\documentclass[10pt]{amsart}

\usepackage[T1]{fontenc}
\usepackage[utf8]{inputenc}
\usepackage{lmodern}
\usepackage{amsmath,amssymb,amsthm,mathtools}
\usepackage{aliascnt}
\usepackage{enumitem}
\usepackage{etoolbox}
\usepackage{needspace}
\usepackage[protrusion=true,expansion=false]{microtype}
\usepackage[colorlinks=true,linkcolor=blue,citecolor=blue,urlcolor=blue]{hyperref}
\usepackage[capitalize,nameinlink]{cleveref}

\allowdisplaybreaks[2]

\makeatletter
\newcommand*\titlelayouterror{\PackageError{title-layout}{Unexpected amsart title format}{}}
\patchcmd{\@maketitle}{\global\topskip42\p@}{\global\topskip\normaltopskip}{}{\titlelayouterror}
\patchcmd{\@setauthors}{\@topsep30\p@}{\@topsep16\p@}{}{\titlelayouterror}
\patchcmd{\@setabstracta}{\skip@20\p@}{\skip@10\p@}{}{\titlelayouterror}
\patchcmd{\@maketitle}{\dimen@34\p@}{\dimen@2\p@}{}{\titlelayouterror}
\makeatother

\newtheorem{theorem}{Theorem}[section]
\newaliascnt{proposition}{theorem}
\newtheorem{proposition}[proposition]{Proposition}
\aliascntresetthe{proposition}
\newaliascnt{lemma}{theorem}
\newtheorem{lemma}[lemma]{Lemma}
\aliascntresetthe{lemma}
\newaliascnt{corollary}{theorem}
\newtheorem{corollary}[corollary]{Corollary}
\aliascntresetthe{corollary}
\theoremstyle{definition}
\newaliascnt{definition}{theorem}

\aliascntresetthe{definition}
\theoremstyle{remark}
\newaliascnt{remark}{theorem}

\aliascntresetthe{remark}

\newcommand{\R}{\mathbb R}

\newcommand{\E}{\mathbb E}
\newcommand{\eps}{\varepsilon}
\newcommand{\one}{\mathbf 1}

\newcommand{\norm}[1]{\lVert #1\rVert}

\title[Linear cotype bounds for scalar expanders]
{Linear Cotype Bounds for Nonlinear Spectral Gaps\\
of Scalar Expanders}

\author{Achintya Raya Polavarapu}
\address{School of Mathematics, Georgia Institute of Technology, Atlanta, GA 30332}
\email{apolavarapu6@gatech.edu}
\date{}

\subjclass[2020]{46B20, 46B42, 46B07, 05C48, 05C50}
\keywords{nonlinear spectral gap, expander graph, cotype, unconditional basis,
Banach lattice, random walk, metric embeddings}

\hypersetup{
  pdftitle={Linear cotype bounds for nonlinear spectral gaps of scalar expanders},
  pdfauthor={Achintya Raya Polavarapu},
  pdfsubject={Linear cotype bounds for scalar expanders and Banach lattices},
  pdfkeywords={nonlinear spectral gap, expander graph, cotype, Banach lattice, random walk}
}

\begin{document}

\begin{abstract}
\setlength{\parindent}{0pt}
\setlength{\parskip}{\medskipamount}
Let \(G\) be a finite simple \(d\)-regular graph, \(d\geq3\), whose second
largest adjacency eigenvalue satisfies \(\lambda_2(G)\leq(1-\delta)d\),
where \(0<\delta\leq1\). We prove that the Poincar\'e constant satisfies
\[
 \gamma(G,\norm{\cdot}_X)\lesssim\frac{qCK}{\delta}
\]
whenever \(X\) is a real Banach space with a \(K\)-unconditional
Schauder basis and cotype \(q\) with constant \(C\). This constant
compares the mean distance over all vertex pairs with the mean distance
across edges. For real Banach lattices the bound is \(O(qC/\delta)\).
The proof estimates where a random walk first leaves each coordinate
level set. As a consequence of the bound, if \(n\)-vertex graphs in a fixed-degree expander
family embed with bounded distortion into spaces with bounded \(C,K\),
then \(q\) must grow at least as fast as \(\log n\). This order is
attained in \(\ell_q^n\).
\end{abstract}
\vspace*{-\dimexpr\headsep+\topskip\relax}
\maketitle

\section{Introduction}

Suppose a single edge joins two equal halves of a graph. Give each
vertex in one half the value zero and each vertex in the other half
the value one. The values differ on half of all ordered pairs of
vertices, but they differ across only one edge.

The scalar Poincar\'e inequality bounds the mean squared difference over
all vertex pairs by a constant times the mean squared difference across
edges. Let \(G=(V,E)\) be a finite connected simple \(d\)-regular graph
on \(n\) vertices. Its normalized spectral gap is
\(\delta(G)=1-\lambda_2(G)/d\), where \(\lambda_2(G)\) is the second
largest adjacency eigenvalue. By the variational formula for this
eigenvalue, the best constant in the scalar inequality is \(1/\delta(G)\)
\cite[Section~1]{NaorComparison}. A family of such graphs is an expander
family if these gaps are bounded away from zero.

What happens to the Poincar\'e constant when the functions take values
in a real normed space \((X,\norm{\cdot}_X)\)?  For \(p\geq1\), the
nonlinear spectral gap
\(\gamma(G,\norm{\cdot}_X^p)\) is the least
\(\gamma\in[0,\infty]\) such that every map \(f:V\to X\) satisfies
\begin{equation}\label{eq:gamma-definition}
 \frac1{n^2}\sum_{u,v\in V}\norm{f(u)-f(v)}_X^p
 \leq
 \frac{\gamma}{|E|}
 \sum_{\{u,v\}\in E}\norm{f(u)-f(v)}_X^p.
\end{equation}
The sum on the left counts ordered pairs; the sum on the right counts
each edge once. We follow \cite{MendelNaor} in writing \(\gamma\)
for this constant. The notation \(\lesssim\) means inequality up to a
universal factor unless a subscript indicates otherwise.

A Schauder basis gives each vector a unique expansion converging in norm.
It is \(K\)-unconditional if changing the signs of the coefficients of
any finite sum increases its norm by at most a factor \(K\). For
\(2\leq q<\infty\), a space \(X\) has cotype \(q\) with constant \(C\)
if, for every finite sequence \((x_s)\) in \(X\) and independent signs
\(\varepsilon_s\) taking the values \(1,-1\) with equal probability,
\begin{equation}\label{eq:cotype}
 \left(\sum_s\norm{x_s}_X^q\right)^{1/q}
 \leq C\left(\E\norm{\sum_s\varepsilon_s x_s}_X^q\right)^{1/q}.
\end{equation}
Finite cotype means that this inequality holds for some finite \(q\) and \(C\).

For spaces with an unconditional basis and finite cotype, a scalar gap
already gave a bound independent of the number of vertices. How does
this bound depend on \(q\)? Even for a \(1\)-unconditional basis, the
earlier bound on \(\gamma(G,\norm{\cdot}_X^2)^{1/2}\) grew
super-exponentially in \(q\) for fixed \(C\) and fixed
gap~\cite[Theorem~A.2]{ADTT}. Under the same graph hypothesis we obtain
a linear bound, both for this square root and for ordinary distances.

\Needspace{16\baselineskip}
\begin{theorem}\label{thm:main}
Let \(G=(V,E)\) be a finite simple \(d\)-regular graph, \(d\geq3\), and suppose that
for some \(0<\delta\leq1\),
\begin{equation}\label{eq:scalar-gap}
 \lambda_2(G)\leq(1-\delta)d.
\end{equation}
If \(X\) is a real Banach space with a \(K\)-unconditional Schauder basis
and cotype \(2\leq q<\infty\) with constant \(C\), then
\begin{equation}\label{eq:main-p1}
 \gamma(G,\norm{\cdot}_X)\lesssim \frac{qCK}{\delta}.
\end{equation}
Equivalently, every map \(f:V\to X\) satisfies
\[
 \frac{1}{|V|^2}\sum_{u,v\in V}\norm{f(u)-f(v)}_X
 \lesssim \frac{qCK}{\delta}\,
 \frac{1}{|E|}\sum_{\{u,v\}\in E}\norm{f(u)-f(v)}_X .
\]
More generally, for \(1\leq p<\infty\),
\begin{equation}\label{eq:main-p}
 \gamma(G,\norm{\cdot}_X^p)^{1/p}\lesssim\frac{p qCK}{\delta}.
\end{equation}
\end{theorem}

A real Banach lattice is an ordered real Banach space in which each
pair of vectors has a supremum and an infimum, and
\(|x|\leq |y|\) implies \(\norm{x}\leq\norm{y}\), where
\(|x|=\sup\{x,-x\}\). The lattice need not have a basis: it suffices
to find a subspace with an unconditional basis containing the finitely
many values of the map.
\begin{corollary}\label{thm:banach-lattice}
If \(G\) satisfies the hypotheses of \Cref{thm:main} and \(X\) is a real
Banach lattice of cotype \(q\) with constant \(C\), then
\begin{equation}\label{eq:lattice-p}
 \gamma(G,\norm{\cdot}_X^p)^{1/p}\lesssim \frac{p qC}{\delta}
 \qquad(1\leq p<\infty).
\end{equation}
\end{corollary}
\begin{proof}
Fix \(f:V\to X\). By Maurey's theorem~\cite[Proposition~1]{Maurey},
in the form used in \cite[Remark~1.9]{ADTT}, the span of \(f(V)\) is
contained in a finite-dimensional subspace \(Y\subseteq X\) with a
\(2\)-unconditional basis. In the inherited norm, \(Y\) has cotype \(q\)
with constant \(C\), and all distances between the values of \(f\) are
the same in \(Y\) as in \(X\). Applying \Cref{thm:main} to \(f\) in
\(Y\) with \(K=2\) gives \eqref{eq:lattice-p}.
\end{proof}

Fix \(d\geq3\) and let \(G\) be a uniformly random simple \(d\)-regular
graph on \(n\) labeled vertices. Friedman's second eigenvalue theorem~\cite{Friedman}
gives \(\lambda_2(G)\leq 2\sqrt{d-1}+\eps\) for every fixed \(\eps>0\)
with probability tending to one as \(n\to\infty\) with \(dn\) even.
Since \(2\sqrt{d-1}<d\), we can choose \(\eps\) so that the normalized
spectral gap is bounded below by a positive constant depending only on
\(d\). On this event \Cref{thm:main} gives
\(\gamma(G,\norm{\cdot}_X)\lesssim_d qCK\) simultaneously for all such \(X\).
This gives the linear dependence on \(q\) for random regular graphs
asked for by Altschuler, Dodos, Tikhomirov, and Tyros~\cite[Problem~6.1]{ADTT}.

\subsection*{The first-exit argument}
To bound a coordinate at a vertex by changes across edges, we follow
a walk from that vertex until the coordinate changes. First subtract a
median in each coordinate, so that at most half the vertices
have positive values and at most half have negative values in that
coordinate. The reduction in \Cref{lem:level-sets} gives a map
\(f:V\to\{-1,0,1\}^k\) with the same property. Fix a coordinate
and let \(S\) be the set on which it equals \(1\), or the set on which it
equals \(-1\). When a walk first leaves \(S\), it crosses an edge on
which this coordinate changes.

The difficulty is that many starting vertices may use the same exit edge.
In \Cref{prop:exit-weights} we retain half the exit probability
from \emph{each} starting vertex in \(S\). At a prescribed oriented edge
and time \(t+1\), the total probability retained is at most
\((1+r)r^t/d\), where \(r=(1+\lambda_2(G)/d)/2\).
To satisfy both requirements, we estimate the loss incurred by imposing
this bound for walks started in an arbitrary subset \(B\subseteq S\).
A bound for \(S\) alone would allow individual vertices to lose all their
weight. The bounds for every \(B\) give precisely the cut inequalities
needed for the max-flow min-cut theorem.

At a fixed exit edge \((x,y)\), we compare all coordinates using the
same walk. By symmetry, the probability of going from \(v\) to \(x\) in
\(t\) steps equals the probability of going from \(x\) to \(v\).
Requiring the walk from \(v\) to stay in a sign set can only decrease
this probability. Thus the unrestricted walk from \(x\) gives a common
distribution on starting vertices for all the sign sets.
With the bound on the retained weights, the cotype estimate in
\Cref{lem:random-positive-vector} gives a factor \(r^{t/q}\) at time
\(t+1\). Summing over time gives \((1-r^{1/q})^{-1}\leq 2q/\delta\).

\subsection*{Relation to earlier estimates}
The question of Pisier and Mendel--Naor asks whether every bounded-degree
expander family satisfies \eqref{eq:gamma-definition} with a bound
independent of the graph for each fixed \(p\) and each Banach space \(X\)
of finite cotype. It remains open (see \cite{Pisier,MendelNaor} and
\cite[Problem~1.3]{ADTT}).

To see why finite cotype matters, give \(V\) its shortest-path distance
\(d_G\). The map \(v\mapsto(d_G(v,w))_{w\in V}\) preserves distances in
\(\ell_\infty^n\) with the maximum norm. By the Maurey--Pisier
theorem~\cite{MaureyPisier}, if \(X\) has no finite cotype, it contains
copies of these spaces with distortion bounded independently of \(n\):
after rescaling each embedding, no distance decreases and none increases
by more than a fixed factor. For a connected graph of degree at most \(d\),
the mean graph distance is \(\gtrsim_d\log n\), while every edge has
length one. Thus even for \(p=1\) the ratio of the two averages in
\eqref{eq:gamma-definition} is unbounded as \(n\) grows. Here and below
\(\log\) denotes the natural logarithm.

For spaces of finite cotype with an unconditional basis, a bound
independent of the number of vertices and the target dimension was
already known. In finite dimension, Odell and Schlumprecht's work
\cite{OdellSchlumprecht} gives a bijection between the unit balls of
\(X\) and Euclidean space that is uniformly continuous in both directions,
with bounds independent of the dimension
\cite[Theorem~A.3]{ADTT}. Building on Ozawa's argument~\cite{Ozawa},
Naor~\cite[Theorem~1.10]{NaorComparison} showed how such a comparison
transfers a Poincar\'e bound from Hilbert space to \(X\).
For a \(1\)-unconditional basis this gives an upper bound of order
\((qC)^{O(q)}\delta^{-4}\) for \(\gamma(G,\norm{\cdot}_X^2)^{1/2}\),
as worked out in \cite[Theorem~A.2]{ADTT}. The exponent of \(q\) in
that estimate is itself proportional to \(q\). Taking \(p=2\) in
\Cref{thm:main} gives \(O(qC/\delta)\) for the same quantity.
For separable infinite-dimensional Banach lattices of finite cotype,
Chaatit~\cite[Corollary~2.11]{Chaatit} obtained a similar comparison
between their unit spheres and the unit sphere of \(\ell_1\).

For finite-dimensional targets, the bound in \cite[Theorem~1.6]{ADTT}
is \(\gamma(G,\norm{\cdot}_X)\lesssim q^{10}C^2K^3\) under the stronger
hypothesis of long-range expansion, with the implicit constant depending
on the degree and expansion parameters.
\Cref{thm:main} gives a linear bound in \(q\) on every regular graph
with a prescribed scalar gap, including expander families for which
long-range expansion fails~\cite[Section~2.3]{AltschulerLRE}.

For the sequence space \(X=\ell_q\), with
\(\norm{x}_q=(\sum_j |x_j|^q)^{1/q}\), and \(p=2\),
Naor~\cite[Theorem~1.5]{NaorComparison}
proved the sharper bound
\(\gamma(G,\norm{\cdot}_{\ell_q}^2)\lesssim q^2/\delta\).
For the full class in \Cref{thm:banach-lattice}, the power
\(\delta^{-2}\) at \(p=2\) is necessary. In \Cref{prop:sharpness} we take the
Cartesian product of several copies of a \(3\)-regular expander and use one copy of
\(\ell_q\) for each factor, adding the norms of the components.
Each edge changes just one component, while
the mean image distance adds over all components. This gives a lower
bound of order \(q/\delta(G)\) on the \(p\)th root of the Poincar\'e
constant. The degree grows with the number of factors.

\section{First exits from a set of vertices}\label{sec:first-exits}

We first consider a single set \(S\) and ask where a walk started
in \(S\) first leaves it. Let \(G\) satisfy \eqref{eq:scalar-gap},
and suppose that \(S\) is nonempty and contains at most half the vertices.
Let \(A_G\) be the adjacency matrix of \(G\). Then \(P=A_G/d\) is the
transition matrix of simple random walk on \(G\).
Its submatrix \(T\) on \(S\) describes the walk killed when it leaves
\(S\): \(T^t(v,x)\) is the probability that it stays in \(S\) for \(t\)
steps and is then at \(x\), starting from \(v\). Write \(\one_B\) for the
indicator of a set \(B\). For \(x\in S\), let \(k(x)\) be the fraction
of its neighbors outside \(S\), the probability of leaving \(S\) in one
step from \(x\). We use counting measure on \(S\) and put
\[
 r=\frac{d+\lambda_2(G)}{2d}.
\]
By \eqref{eq:scalar-gap}, \(r<1\). The adjacency eigenvalues other than
the largest sum to \(-d\), so \(\lambda_2(G)\geq-d/(n-1)\) and \(r>0\).

\subsection{A boundary estimate for the killed walk}

Few walks remain in \(S\) after many steps, but this alone does not say
where they leave. To estimate the exits we need to know how much of the
walk reaches vertices with a neighbor outside \(S\).

\begin{lemma}\label{lem:boundary-dissipation}
For every \(h\in\ell_2(S)\), we have \(\norm{Th}_2\leq r\norm{h}_2\) and
\begin{equation}\label{eq:boundary-dissipation}
 \sum_{t\geq0}r^{-t}\sum_{x\in S}k(x)\big(T^th(x)\big)^2
 \leq (1+r)\norm{h}_2^2.
\end{equation}
\end{lemma}

\begin{proof}
Extend \(h\) by zero to \(V\). Its orthogonal projection onto the constants
has squared norm \((\sum_{x\in S}h(x))^2/n\leq\norm{h}_2^2/2\).
On the orthogonal complement of the constants, the Rayleigh quotient of
\(A_G/d\) is at most \(\lambda_2(G)/d\). Thus
\[
 \langle Th,h\rangle\leq
 \frac{\lambda_2(G)}{d}\norm{h}_2^2+
 \left(1-\frac{\lambda_2(G)}{d}\right)
 \frac{\big|\sum_{x\in S}h(x)\big|^2}{n}
 \leq r\norm{h}_2^2.
\]
Applying this to \(|h|\) and
using \(|\langle Th,h\rangle|\leq\langle T|h|,|h|\rangle\) gives the norm bound.

Since the row sums of \(T\) are at most one, \((Th(x))^2\leq T(h^2)(x)\) by
Cauchy--Schwarz. Summing over \(x\) and using symmetry yields
\[
 \sum_x k(x)h(x)^2\leq\norm{h}_2^2-\norm{Th}_2^2.
\]
Apply this inequality to \(T^th\) and put \(u_t=\norm{T^th}_2^2\).
Since \(u_t\leq r^{2t}u_0\), summation by parts gives
\[
 \sum_{t=0}^{N}r^{-t}(u_t-u_{t+1})
 \leq u_0+(1-r)\sum_{t=1}^{N}r^{-t}u_t
 \leq (1+r)u_0.
\]
Letting \(N\to\infty\) proves \eqref{eq:boundary-dissipation}.
\end{proof}

\subsection{Assigning weights to exit edges and times}

In this subsection \(e=(x,y)\) ranges over edges directed from \(S\) to
its complement. A walk starting at \(v\in S\) first leaves through \(e\)
at step \(t+1\) with probability \(T^t(v,x)/d\): it must stay in \(S\)
for \(t\) steps, reach \(x\), and then cross \(e\). Many starts may use
the same edge at the same time. We discard some of these exits while
keeping half of the exit probability from every starting vertex.

\begin{proposition}\label{prop:exit-weights}
Choose an integer \(H\geq1\) with \(\sqrt n\,r^H\leq1/4\).
For \(v\in S\), each edge \(e\) leaving \(S\), and \(0\leq t<H\), there are
weights \(0\leq\theta_v(e,t)\leq T^t(v,x)/d\) such that
\begin{align}
 \sum_{e,t}\theta_v(e,t)&=\tfrac12
       && \text{for each }v\in S,\label{eq:first-exit-coverage}\\
 \sum_{v\in S}\theta_v(e,t)&\leq (1+r)r^t/d
       && \text{for each }(e,t).\label{eq:edge-exit-bound}
\end{align}
\end{proposition}

\begin{proof}
Join a source to each \(v\in S\) with capacity \(1/2\), each \(v\) to
\((e,t)\) with capacity \(T^t(v,x)/d\), and each \((e,t)\) to a sink
with capacity \((1+r)r^t/d\). Fix the set \(B\subseteq S\) of starting
vertices on the source side of a cut. Its capacity is the sum of the
capacities of arcs directed from the source side to the sink side.
Put \(h_t=T^t\one_B\). By symmetry,
the incoming capacity at \((e,t)\) from \(B\) is \(h_t(x)/d\). Placing
\((e,t)\) on the sink side costs this amount. Placing it on the source
side costs its capacity to the sink. We may choose the cheaper side
independently for each \((e,t)\). Each pair therefore contributes the
smaller of its incoming exit mass and \((1+r)r^t/d\).

The probability of surviving \(H\) steps from \(v\) is
\((T^H\one_S)(v)\leq r^H\sqrt{|S|}\leq1/4\). Hence the total exit mass
by time \(H\) from starts in \(B\) is at least \(3|B|/4\).
There are \(dk(x)\) exit edges from \(x\). Writing \(u_+=\max\{u,0\}\),
we sum the excess over these edges:
\[
 \sum_{t=0}^{H-1}\sum_{x\in S} k(x)\big(h_t(x)-(1+r)r^t\big)_+
 \leq \frac{1}{4(1+r)}\sum_{t\geq0}r^{-t}\sum_{x\in S}k(x)h_t(x)^2
 \leq \frac{|B|}{4}.
\]
The first inequality uses \((u-a)_+\leq u^2/(4a)\) for \(u\geq0\), \(a>0\),
and the second uses \eqref{eq:boundary-dissipation} with \(h=\one_B\).
The pairs \((e,t)\) therefore contribute at least \(|B|/2\). The source arcs
to \(S\setminus B\) contribute \(|S\setminus B|/2\), so every cut has
capacity at least \(|S|/2\). By the max-flow min-cut theorem~\cite{FordFulkerson},
the maximum flow equals the minimum cut capacity. The source arcs have
total capacity \(|S|/2\), so a maximum flow saturates each one. Its value from \(v\) to
\((e,t)\) gives \(\theta_v(e,t)\).
\end{proof}

\section{Cotype and coordinate level sets}\label{sec:setup}

The exit estimate concerns one set at a time. To use it for a map into
\(X\), we take level sets in each coordinate and use cotype to add the
weights from the different starting vertices.

\subsection{Cotype on positive vectors}
For a \(1\)-unconditional norm \(U\) on \(\R^k\), we have
\(U(x)=U(|x|)\) and \(U(x)\leq U(y)\) whenever \(0\leq x\leq y\).
Here \(|x|\) is the vector of absolute values, and all vector inequalities
are coordinatewise. Since exit weights are nonnegative, the only
consequence of cotype we need is
\begin{equation}\label{eq:positive-estimate}
 \left(\sum_{s=1}^m U(x_s)^q\right)^{1/q}
 \leq C\,U\left(\sum_{s=1}^m x_s\right)
 \qquad(x_1,\ldots,x_m\in\R_+^k).
\end{equation}
For disjointly supported vectors this is the usual lower \(q\)-estimate.
We can make the norm invariant under sign changes and still use \(C\)
in \eqref{eq:positive-estimate}, even when the supports overlap. The
factor \(K\) enters only when we return to the original norm.

\begin{lemma}\label{lem:positive-renorming}
Suppose the coordinate basis of \((\R^k,\norm{\cdot}_X)\) is
\(K\)-unconditional and \(X\) has cotype \(q\) with constant \(C\). Then
\[
 U(x)=\max\{\norm{y}_X: |y|=|x|\}
\]
is a \(1\)-unconditional norm satisfying \eqref{eq:positive-estimate},
and \(\norm{x}_X\leq U(x)\leq K\norm{x}_X\).
\end{lemma}

\begin{proof}
The maximum is over all sign changes of \(x\), so \(U\) is a norm invariant
under sign changes. The norm comparison follows from \(K\)-unconditionality.
For each \(x_s\geq0\),
choose a signed vector \(y_s\) with \(|y_s|=x_s\) and \(\norm{y_s}_X=U(x_s)\).
For every choice of signs \(\varepsilon_s\),
\(\big|\sum_s\varepsilon_s y_s\big|\leq\sum_s x_s\). Thus
\[
 \left(\sum_s U(x_s)^q\right)^{1/q}
 \leq C\left(\E\norm{\sum_s\varepsilon_s y_s}_X^q\right)^{1/q}
 \leq C\,U\left(\sum_s x_s\right),
\]
where the last inequality uses the monotonicity of \(U\).
\end{proof}

\subsection{Positive vectors with overlapping supports}

Several starting vertices can contribute to the same coordinate at an
exit edge. \Cref{prop:exit-weights} bounds both the weight from each
start and the sum of these weights. After normalization these give
\(h\leq z\) and \(\E h\leq az\) in the following lemma. The weights
play the role of the support counts in \cite[Proposition~3.4]{ADTT}.
We need the bound to decrease with \(a\), since later exit times have
smaller total weight.

\begin{lemma}\label{lem:random-positive-vector}
Let \(U\) be a \(1\)-unconditional norm satisfying
\eqref{eq:positive-estimate}. Suppose that \(h\) is a random vector in
\(\R_+^k\), \(0\leq h\leq z\) almost surely, and \(\E h\leq a z\)
for some \(z\in\R_+^k\) and \(a>0\). Then
\[
 \bigl(\E U(h)^q\bigr)^{1/q}\leq 2Ca^{1/q}U(z).
\]
\end{lemma}

\begin{proof}
Take \(N=\lceil 1/a\rceil\) independent copies \(h_1,\ldots,h_N\)
of \(h\). For each \(i\), the other \(N-1\) copies together contribute
at most \((N-1)a z_j\leq z_j\) in expectation to coordinate \(j\).
We divide the available amount \(z_j\) among the copies so that
each copy keeps, in expectation, at least half its original weight.
Set \(Y_{i,j}=0\) if \(z_j=0\). For \(z_j>0\), put
\[
 Y_{i,j}=\frac{z_j h_{i,j}}{z_j+\sum_{s\ne i}h_{s,j}}.
\]
As \(h_{i,j}\leq z_j\), the denominator is at least
\(\sum_s h_{s,j}\). If that sum is positive, summing over \(i\) gives
\(\sum_iY_{i,j}\leq z_j\). If it is zero, each \(Y_{i,j}\) is zero.
Thus \(\sum_iY_i\leq z\), and \eqref{eq:positive-estimate} gives
\(\sum_i U(Y_i)^q\leq C^q U(z)^q\) even when the supports overlap.

Condition on \(h_i\). For \(z_j>0\), independence and convexity of
the reciprocal on the positive reals give
\[
 \E[Y_{i,j}\mid h_i]
 \geq\frac{z_j h_{i,j}}{z_j+(N-1)(\E h)_j}
 \geq\frac{h_{i,j}}{1+(N-1)a}\geq\tfrac12 h_{i,j}.
\]
Monotonicity of \(U\) and conditional Jensen therefore give
\(\E U(Y_i)^q\geq 2^{-q}\E U(h)^q\).
Summing over \(i\) yields \(N2^{-q}\E U(h)^q\leq C^qU(z)^q\).
Now use \(1/N\leq a\).
\end{proof}

\subsection{Reduction to signed indicators}
We want to replace each real coordinate by signed indicators without
changing the distance between any two values. The sign sets for these
indicators must also be small enough for the exit estimate. Subtracting
a coordinatewise median gives a map \(g\) satisfying
\begin{equation}\label{eq:median-condition}
 |\{v:g(v)_j>0\}|\leq n/2,\qquad |\{v:g(v)_j<0\}|\leq n/2
 \quad(1\leq j\leq k).
\end{equation}
Following \cite[Corollary~3.3]{ADTT}, we reduce to maps with coordinates
in \(\{-1,0,1\}\). We use the values of \(|g(v)_j|\) as thresholds
and weight each by its distance from the preceding one.

\begin{lemma}\label{lem:level-sets}
Let \(g:V\to\R^k\) and let \(U\) be a \(1\)-unconditional norm satisfying
\eqref{eq:positive-estimate}. There is a map \(F\) from \(V\) into a
finite-dimensional coordinate space, with every coordinate taking values
in \(\{-1,0,1\}\). This space admits a \(1\)-unconditional norm
\(\widetilde U\) satisfying \eqref{eq:positive-estimate} with the same
\(q,C\), for which
\begin{equation}\label{eq:level-distances}
 \widetilde U(F(v))=U(g(v)),\qquad
 \widetilde U(F(u)-F(v))=U(g(u)-g(v))\quad(u,v\in V).
\end{equation}
If \(g\) satisfies \eqref{eq:median-condition}, then so does \(F\).
\end{lemma}

\begin{proof}
The assertion is immediate for the zero map. Otherwise, list the distinct
positive values of \(|g(v)_j|\) over all vertices and coordinates as
\(t_1<\cdots<t_m\). Put \(t_0=0\) and \(w_s=t_s-t_{s-1}\).
Use these same thresholds in every coordinate and define
\[
 F(v)_{j,s}=\operatorname{sgn}(g(v)_j)\one_{\{|g(v)_j|\geq t_s\}},
 \qquad
 \widetilde U(y)=U\left(\left(\sum_s w_s|y_{j,s}|\right)_{j=1}^k\right).
\]
The positive weights and monotonicity of \(U\) make \(\widetilde U\) a
\(1\)-unconditional norm on \(\R^{km}\). For nonnegative vectors the
weighted sums in each coordinate are additive. Applying
\eqref{eq:positive-estimate} to the vectors of weighted sums therefore
gives the same estimate for \(\widetilde U\), with the same \(q,C\).

For fixed \(j\), the sets of thresholds below two absolute values are
nested. If the two coordinate values have the same sign, the differing
indicators have total weight equal to the difference of their absolute
values. For opposite signs, each crossed level contributes once for each
endpoint that reaches it. In either case,
\[
 \sum_s w_s|F(u)_{j,s}-F(v)_{j,s}|=|g(u)_j-g(v)_j|.
\]
Comparison with zero gives \(\sum_s w_s|F(v)_{j,s}|=|g(v)_j|\),
proving \eqref{eq:level-distances}. Wherever a threshold coordinate is
nonzero, it has the sign of the original coordinate. Thus neither sign
occurs more often in a threshold coordinate than in the original one,
which proves the assertion about medians.
\end{proof}

It now suffices to bound distances to the median for maps taking values
in \(\{-1,0,1\}^k\). The triangle inequality then gives a bound for
all vertex pairs, and the same estimate also controls powers of the distance.
\begin{proposition}\label{prop:reduction}
Fix \(q\geq2\), \(C\geq1\), and a \(d\)-regular graph \(G\). Suppose
there is a constant \(\Pi\) such that in every dimension, for every \(1\)-unconditional norm \(U\)
satisfying \eqref{eq:positive-estimate}, one has
\begin{equation}\label{eq:median-poincare}
 \sum_{v\in V}U(f(v))\leq\Pi\sum_{\{u,v\}\in E}U(f(u)-f(v))
\end{equation}
for every \(f:V\to\{-1,0,1\}^k\) satisfying \eqref{eq:median-condition}.
Then \eqref{eq:median-poincare} holds for every \(f:V\to\R^k\)
satisfying \eqref{eq:median-condition}, and
\begin{equation}\label{eq:extrapolation}
 \gamma(G,U)\leq d\Pi,\qquad
 \gamma(G,U^p)\leq\tfrac12(2pd\Pi)^p\qquad(p\geq1).
\end{equation}
\end{proposition}

\begin{proof}
Apply the hypothesis to the map and norm in \Cref{lem:level-sets}.
The identities \eqref{eq:level-distances} give the first assertion.
For an arbitrary map, subtract a coordinatewise median and denote the
result by \(g\). For \(p=1\), summing the triangle inequality over all
ordered pairs gives \(\sum_{u,v}U(g(u)-g(v))\leq2n\sum_v U(g(v))\).
As \(|E|=dn/2\), this proves \(\gamma(G,U)\leq d\Pi\).

For \(p>1\) we use Matou\v{s}ek's extrapolation \cite{Matousek} in
the form of \cite[proof of Corollary~3.3]{ADTT}. Put
\(T_p(x)=U(x)^{p-1}x\). The map \(T_p\circ g\) still satisfies
\eqref{eq:median-condition}, and
\begin{equation}\label{eq:radial-extrapolation}
 U(T_p(x)-T_p(y))\leq p\max\{U(x),U(y)\}^{p-1}U(x-y).
\end{equation}
Indeed, if \(a=U(x)\geq b=U(y)>0\), then
\[
 U(T_p(x)-T_p(y))\leq a^{p-1}U(x-y)+b(a^{p-1}-b^{p-1}).
\]
On \(b\leq t\leq a\) one has \(bt^{p-2}\leq t^{p-1}\leq a^{p-1}\).
Integrating gives \(b(a^{p-1}-b^{p-1})\leq(p-1)a^{p-1}(a-b)\).
Since \(a-b\leq U(x-y)\),
\eqref{eq:radial-extrapolation} follows. The case \(b=0\) is immediate.

Write \(S=\sum_v U(g(v))^p\). Apply \eqref{eq:median-poincare} to
\(T_p\circ g\), use \eqref{eq:radial-extrapolation}, and apply H\"older
on the edges. This gives
\[
 S\leq p\Pi\,(dS)^{1-1/p}
 \left(\sum_{\{u,v\}\in E}U(g(u)-g(v))^p\right)^{1/p},
\]
since \(\sum_{\{u,v\}\in E}\max\{U(g(u)),U(g(v))\}^p\leq dS\).
For \(S>0\), divide by \(S^{1-1/p}\) and raise to the power \(p\).
The bound \(\sum_{u,v}U(g(u)-g(v))^p\leq2^p n S\) and the identity
\(|E|=dn/2\) now give \eqref{eq:extrapolation}. The case \(S=0\) is trivial.
\end{proof}

\section{The Poincar\'e inequality}\label{sec:analytic}

We count the exit weights in two ways. At a starting vertex they account
for half of each nonzero coordinate in absolute value. At an exit edge
we must bound the sum of their norms by how much the map changes across
that edge.

\subsection{Walks of fixed length ending at a given oriented edge}

Let \(G\) satisfy \eqref{eq:scalar-gap}. Fix a \(1\)-unconditional norm \(U\) satisfying \eqref{eq:positive-estimate}
and a map \(f:V\to\{-1,0,1\}^k\) with median zero in each coordinate.
Apply \Cref{prop:exit-weights} to the positive and negative sets of each
coordinate, omitting empty sets. Each contains at most \(n/2\) vertices,
so the same \(H\) with \(\sqrt n\,r^H\leq1/4\) works for all of them.
Let \(a_v(e,t)\) be the vector
of exit weights from \(v\) through the oriented edge \(e\) at time \(t+1\).
For coordinate \(j\), use the weight for the sign set containing \(v\).
The entry is zero if \(f(v)_j=0\) or if \(e\) does not leave that set.
Summing over times and oriented edges, \eqref{eq:first-exit-coverage} gives
\begin{equation}\label{eq:sum-at-start}
 \sum_{t=0}^{H-1}\sum_e a_v(e,t)=\tfrac12|f(v)|.
\end{equation}

Fix \(e=(x,y)\) and \(t\), and write \(a_v=a_v(e,t)\). For the
unrestricted walk put \(\pi(v)=P^t(x,v)\). By symmetry,
\(\pi(v)=P^t(v,x)\). Thus the same probability distribution on \(V\) can
be used for all coordinates, even though their sign sets differ. Put
\(z=|f(x)-f(y)|\). For each \(j\), only the sign set containing \(x\)
can contribute, and then \(y\) lies outside it and \(z_j\geq1\).
For this set and \(v\) in it, the killed walk \(T\) gives
\[
 d(a_v)_j\leq T^t(v,x)\leq P^t(v,x)=\pi(v).
\]
The first inequality is the bound on the retained exit probability.
The second holds because every walk counted by \(T^t(v,x)\) is also
counted by \(P^t(v,x)\). Together with \eqref{eq:edge-exit-bound},
this gives the coordinatewise bounds
\[
 d a_v\leq \pi(v)z,\qquad
 d\sum_v a_v\leq(1+r)r^t z.
\]
If \(\pi(v)=0\), the first bound forces \(a_v=0\). For the remaining
vertices take the vector \(h=d a_v/\pi(v)\) with probability \(\pi(v)\).
Then \(0\leq h\leq z\) and \(\E h=d\sum_v a_v\leq(1+r)r^t z\).

Using \(\E U(h)\leq(\E U(h)^q)^{1/q}\) and \Cref{lem:random-positive-vector},
we obtain
\begin{equation}\label{eq:sum-at-edge}
 \sum_v U(a_v)=\frac1d\E U(h)
 \leq\frac{2C}{d}(1+r)^{1/q}r^{t/q}U(f(x)-f(y)).
\end{equation}

Take norms in \eqref{eq:sum-at-start}, apply the triangle
inequality, and sum over \(v\). Using \eqref{eq:sum-at-edge} for each
\(t\) and both orientations of every edge, we obtain
\begin{equation}\label{eq:indicator-poincare}
 \sum_{v\in V}U(f(v))
 \leq\frac{8C}{d}\frac{(1+r)^{1/q}}{1-r^{1/q}}
 \sum_{\{x,y\}\in E}U(f(x)-f(y)).
\end{equation}
The denominator comes from summing \(r^{t/q}\) over \(t\geq0\). In particular
the bound is independent of the cutoff \(H\).

\subsection{Passage to arbitrary maps}

The sum over exit times is where \(q\) enters: the weights decay as
\(r^{t/q}\) after applying cotype. It remains to estimate this sum and
pass from signed indicators to the original map.

\begin{proof}[Proof of \Cref{thm:main}]
In finite dimension take \(U\) from \Cref{lem:positive-renorming}. By
\eqref{eq:indicator-poincare}, \Cref{prop:reduction} applies with
\[
 d\Pi=8C\frac{(1+r)^{1/q}}{1-r^{1/q}}.
\]
Concavity gives \(1-r^{1/q}\geq(1-r)/q\), while \eqref{eq:scalar-gap}
gives \(1-r\geq\delta/2\). Since \(q\geq2\), we obtain \(d\Pi\leq32qC/\delta\).
Now \(\norm{x}_X\leq U(x)\leq K\norm{x}_X\), so \eqref{eq:extrapolation}
yields
\begin{align*}
 \gamma(G,\norm{\cdot}_X)&\leq dK\Pi\leq 32qCK/\delta,\\
 \gamma(G,\norm{\cdot}_X^p)&\leq\tfrac12(2pdK\Pi)^p
 \leq\tfrac12(64pqCK/\delta)^p.
\end{align*}
For a Schauder basis, the span of the first \(N\) basis vectors inherits
both the cotype and the unconditionality bounds. Apply the finite-dimensional
inequalities to the corresponding partial sums of \(f\) at each vertex, and
let \(N\to\infty\). Both sides converge since \(V\) is finite.
\end{proof}

\section{Sharpness and distortion}\label{sec:applications}

For \(2\leq q\lesssim\log n\), Matou\v{s}ek's embedding
theorem~\cite{Matousek} gives an embedding of the shortest-path metric
of a connected \(n\)-vertex graph into \(\ell_q\) with distortion
\(O((\log n)/q)\). At the scale where the least ratio of image distance
to graph distance is one, the mean image distance for a \(d\)-regular graph is
\(\gtrsim_d\log n\), while every image edge has length \(O((\log n)/q)\).
Thus \(\gamma(G,\norm{\cdot}_{\ell_q})\gtrsim_d q\), as observed in
\cite[Remark~1.7]{ADTT}. To see the dependence on the scalar gap as well,
we take a product of several copies of the same graph. Adding copies
reduces the gap. Each edge still changes only one coordinate, while
the distances between typical pairs add over all coordinates.

\begin{proposition}\label{prop:sharpness}
There are cubic expanders \(H\) of arbitrarily large order \(N\) whose
normalized gaps are at most one and bounded below by an absolute positive constant.
Put \(q=\log N\). For each integer \(L\geq1\), let \(G\) be the Cartesian
product of \(L\) copies of \(H\). Its vertices are \(L\)-tuples of vertices
of \(H\). Two tuples are adjacent if they agree in all but one coordinate
and their entries in that coordinate are adjacent in \(H\).
Let \(X\) be the space of \(L\)-tuples of vectors in \(\ell_q^N\),
with norm the sum of the \(\ell_q\) norms of the entries. Then
\(\delta(G)=\delta(H)/L\), and the
Banach lattice \(X\) has a \(1\)-unconditional basis and cotype \(q\)
with constant at most \(\sqrt2\). Moreover,
\[
 \gamma(G,\norm{\cdot}_X^p)^{1/p}\gtrsim\frac{q}{\delta(G)}
 \qquad(1\leq p<\infty).
\]
\end{proposition}

\begin{proof}
Taking \(d=3\) in Friedman's theorem gives simple cubic graphs \(H\) of arbitrarily
large order \(N\) with normalized gaps bounded away from zero. Take \(N\) large
enough that \(q=\log N\geq2\) and that \(H\) has two nonadjacent vertices.
The difference of their coordinate vectors has mean zero and adjacency
Rayleigh quotient zero, so \(\delta(H)\leq1\).

Define \(f:V(H)\to\ell_q^N\)
by \(f(u)=(d_H(u,w))_{w\in V(H)}\), where \(d_H\) is shortest-path distance in \(H\).
Each coordinate is
\(1\)-Lipschitz, and the coordinate \(w=u\) realizes \(d_H(u,v)\). Hence
\[
 d_H(u,v)\leq\norm{f(u)-f(v)}_q\leq N^{1/q}d_H(u,v)=e\,d_H(u,v).
\]
A ball of radius \(R\) in a cubic graph has at most \(1+3(2^R-1)\)
vertices. For \(R=\lfloor\frac12\log_2 N\rfloor\), this is at most \(N/2\)
for all sufficiently large \(N\). At least half the vertices therefore
lie beyond distance \(R\) from any given vertex. The mean graph distance,
and hence the mean distance between the \(f(u)\), is \(\gtrsim\log N=q\).

Define \(F:V(G)\to X\) by
\[
 F(u_1,\ldots,u_L)=(f(u_1),\ldots,f(u_L)).
\]
There are three choices of neighbor in each coordinate, so \(G\) is
\(3L\)-regular. Its adjacency eigenvalues are sums of \(L\)
eigenvalues of \(H\), so its second largest eigenvalue is
\(3(L-1)+\lambda_2(H)\) and \(\delta(G)=\delta(H)/L\).
An edge of \(G\) changes one component of \(F\), and hence has image length
at most \(e\). For independent uniform vertices \(u,v\), the mean of
\(\norm{F(u)-F(v)}_X\) is \(L\) times the corresponding mean in one copy.
By Jensen's inequality and \eqref{eq:gamma-definition},
\[
 \gamma(G,\norm{\cdot}_X^p)^{1/p}
 \geq \frac{\E\norm{F(u)-F(v)}_X}{e}
 \gtrsim Lq
 =\frac{\delta(H)q}{\delta(G)}
 \gtrsim \frac{q}{\delta(G)}.
\]

It remains to check the cotype constant of \(X\) uniformly in \(L\) and \(q\).
For \(y_s\in\ell_q^N\), use \(\norm{\cdot}_q\leq\norm{\cdot}_2\) on the index
\(s\), then the scalar Khintchine inequality~\cite{Szarek} in each coordinate,
and finally Minkowski's inequality. This gives
\begin{align*}
 \left(\sum_s\norm{y_s}_q^q\right)^{1/q}
 &\leq \left\|\left(\sum_s |y_s|^2\right)^{1/2}\right\|_q
 \leq \sqrt2\left\|\E\left|\sum_s\varepsilon_s y_s\right|\right\|_q
 \leq\sqrt2\,\E\norm{\sum_s\varepsilon_s y_s}_q .
\end{align*}
For \(x_s=(x_{s,a})_{a=1}^L\in X\), Minkowski's inequality on the index \(s\)
now gives
\[
 \left(\sum_s\norm{x_s}_X^q\right)^{1/q}
 \leq\sum_{a=1}^L\left(\sum_s\norm{x_{s,a}}_q^q\right)^{1/q}
 \leq\sqrt2\,\E\norm{\sum_s\varepsilon_s x_s}_X .
\]
This implies \eqref{eq:cotype} with constant at most \(\sqrt2\). The
standard coordinate basis is \(1\)-unconditional.
\end{proof}

Each of these finite-dimensional targets also has cotype \(2\), but the
corresponding constants grow with \(N\). Indeed, applying \eqref{eq:cotype}
to the \(N\) coordinate vectors in one \(\ell_q^N\) summand gives a
cotype-\(2\) constant of at least \(\sqrt N/N^{1/q}=\sqrt N/e\).

In a graph of bounded degree, the mean distance between vertices grows
with \(n\), while every edge has length one. If an embedding does not
shorten distances, the Poincar\'e inequality forces it to stretch some edge.
\begin{corollary}\label{cor:distortion}
Under the hypotheses of \Cref{thm:main}, every bi-Lipschitz embedding
of \((V,d_G)\) into \(X\) has distortion
\[
 D\gtrsim_d\frac{\delta\log n}{qCK}.
\]
For a Banach lattice of cotype \(q\) with constant \(C\), the denominator
can be replaced by \(qC\).
\end{corollary}
\begin{proof}
Scale the map so that \(d_G(u,v)\leq\norm{f(u)-f(v)}_X\leq Dd_G(u,v)\).
The ball count used above, with
\(1+d\sum_{i=0}^{R-1}(d-1)^i\) in place of \(1+3(2^R-1)\), gives mean graph
distance \(\gtrsim_d\log n\). The mean image distance is therefore
\(\gtrsim_d\log n\), while each image edge has length at most \(D\).
Apply \eqref{eq:main-p1}, or \eqref{eq:lattice-p} with \(p=1\).
\end{proof}

For a fixed-degree expander family, consider embeddings of distortion
less than \(20\) into finite-dimensional spaces with \(C,K\leq20\).
Then \Cref{cor:distortion} gives \(q\gtrsim_{d,\delta}\log n\).
This order is attained. For \(q=\log n\geq2\), the distance map
\(v\mapsto(d_G(v,w))_{w\in V}\) into \(\ell_q^n\) has distortion at
most \(n^{1/q}=e\). The space \(\ell_q^n\) has a \(1\)-unconditional basis
and cotype \(q\) with constant at most \(\sqrt2\), as in \Cref{prop:sharpness}.
Thus the infimum of the admissible cotype exponents has order \(\log n\).
For random regular graphs of fixed degree \(d\geq6\), the earlier
lower bound on this quantity in expectation was of order
\((\log n)^{1/10}\)~\cite[Corollary~1.13]{ADTT}.

\section*{Acknowledgment}

In the original draft, the author's proof of the main theorem distributed
the weight at each exit edge and time among individual walks. During the
revision, the author used ChatGPT Pro and OpenAI Codex to explore other
arguments, check details of the proofs, and revise the exposition.
Codex suggested comparing the exit probabilities for all coordinates with
the transition probabilities of a single unrestricted random walk. This
comparison permits a direct application of the cotype inequality to the
exit weights and led to a shorter proof and a shorter manuscript.
The author takes responsibility for the arguments and for any errors
that remain.

\bibliographystyle{alpha}
\bibliography{bibliography}

\end{document}